\documentclass[10pt,a4paper]{article}
\usepackage[latin1]{inputenc}
\usepackage{amsmath}
\usepackage{amsfonts,amssymb}
\usepackage{epsfig}
\usepackage{psfrag}
\usepackage{mathrsfs}

\def\dsp{\displaystyle}

\def\R{\mathbb{R}}

\def\T{\mathbb{T}}

\def\b0{\boldsymbol{0}}
\def\bA{\boldsymbol{A}}
\def\bB{\boldsymbol{B}}

\def\bF{\boldsymbol{F}}
\def\bbf{\boldsymbol{f}}

\def\bG{\boldsymbol{G}}

\def\bu{\boldsymbol{u}}
\def\bv{\boldsymbol{v}}
\def\bx{\boldsymbol{x}}

\def\bdelta{\boldsymbol{\delta}}

\def\btau{\boldsymbol{\tau}}

\def\div{\mbox{\,{\textrm{div}}}}
\def\transp{\,^T\!}
\newcommand{\proof}{{\bf Proof~~}}
\newcommand{\cqfd}{{\nobreak\hfil\penalty50\hskip2em\hbox{}\nobreak\hfil $\square$\qquad\parfillskip=0pt\finalhyphendemerits=0\par\medskip}}

\newtheorem{theorem}{Theorem}[section]
\newtheorem{lemma}{Lemma}[section]
\newtheorem{definition}{Definition}[section]
\newtheorem{remark}{Remark}[section]

\newtheorem{proposition}{Proposition}[section]

\begin{document}

\setlength\parindent{0pt}

\begin{center}
{ \Large \bf Global existence results\\[0.0cm] for some viscoelastic models\\[0.25cm] with an integral constitutive law}
\end{center}

\begin{center}
Laurent Chupin
\footnote{
Laboratoire de Math\'ematiques, UMR 6620 \\
Universit\'e Blaise Pascal, Campus des C\'ezeaux, F-63177 Aubi\`ere Cedex, France.\\
{\it laurent.chupin@math.univ-bpclermont.fr}
}
\end{center}

\begin{abstract}
We provide a proof of global regularity of solutions of some models of viscoelastic flow with an integral constitutive law, in the two spatial dimensions and in a periodic domain.
Models that are included in these results are classical models for flow memory: for instance some K-BKZ models, the PSM model or the Wagner model. The proof is based on the fact that these models naturally give a $L^\infty$-bound on the stress and that they allow to control the spatial gradient of the stress. The main result does not cover the case of the Oldroyd-B model.
\end{abstract}



\section{Introduction}\label{part:introduction}

\subsection{Presentation of the result}

In this article, we are interested in the global (with respect to the time) existence result for models of viscoelastic fluids. Usually, obtain a global existence result for highly nonlinear system of PDE is quite challenging. The models we are interested in here are nonlinear at several levels:
The first one is the well-known nonlinearity of the Navier-Stokes equations describing the hydrodynamics, this is the main reason we do not expect to have results in the three dimensional case.
The second level of nonlinearity comes from to the rheology that we consider.
More precisely, the viscoelasticity is described by the constitutive relation linking the stress and the strain.
The framework for our study corresponds to the case where the extra-stress~$\btau$ is given, at any time~$t$ and at each point~$\bx$ by an integral law of the form
\begin{equation}\label{stress-integral}
\btau(t,\bx) = \int_{-\infty}^t \mathcal F(t-\sigma,\bF(\sigma,t,\bx))\, \mathrm d\sigma.
\end{equation}
The tensor $\bF$ contains all the information of past deformations.
It naturally depends on the velocity field of the flow: the relation~\eqref{stress-integral} is then strongly coupled with the Navier-Stokes equations. 
Under assumptions on the behavior of the functional $\mathcal F$, we prove that the resulting system admits a global solution, in the two-dimensional case and in a periodic domain, but without assuming that the data are small.

These assumptions on the functional~$\mathcal F$ (see the Subsection~\ref{sec:assumption} and the Remark~\ref{rem:non-separable-stress}) allow us to consider most of the usual integral models: the Wagner model, the PSM model and some K-BKZ models.

\subsection{Mathematical results on viscoelastic model with an integral constitutive law}

The integral models are extensively studied in the last half a century.
In this regard, we can read the review article written by Mitsoulis~\cite{Mitsoulis} for the $50$th anniversary of the K-BKZ models.
However, there are few mathematical works on such viscoelastic models.
The first significant results are probably due to Kim~\cite{Kim}, M. Renardy~\cite{Renardy3}, Hrusa and Renardy~\cite{Renardy2}, Hrusa, Nohel and Renardy~\cite[Section IV.5]{Renardy1}.
Kim discusses a situation in which the nonlinearity in the constitutive equation (that is the functional~$\mathcal F$ in the relation~\eqref{stress-integral}) has special form. Renardy, Hrusa and Nohel study spatially periodic three dimensional motions with a more general nonlinearity (but sufficiently smooth).
In all these works the solution is either local in time or global but with small data.
Later, Brandon and Hrusa~\cite{Brandon-Hrusa} study a one dimensional model with a singularity in the nonlinearity: they obtain global existence results for sufficiently small data.
Very recently - see~\cite{Chupin3} - some theoretical results are proved for a large family of nonlinearities: local existence, global existence with small data and uniqueness results.

\subsection{Some global existence results for viscoelastic models}

\paragraph{The Oldroyd-B model}

There are many ways to describe a flow of viscoelastic fluid. The most famous model is the Oldroyd-B model for which the question of global existence is still open, even in the two dimensional case.
This model expresses the constitutive relation between the extra-stress~$\btau$ and the deformation tensor $D\bu=\frac{1}{2}(\nabla \bu+ \transp{(\nabla \bu)})$ as follows:
$$\lambda \stackrel{\triangledown}{\btau} + \btau = 2 \mu D \bu.$$
In this expression, the constants~$\mu$ and~$\lambda$ respectively correspond to a polymeric viscosity and a relaxation time. The notation {\footnotesize $\triangledown$} stands for the upper-convective derivative:
$$\stackrel{\triangledown}{\btau} \, = \partial_t\btau + \bu\cdot \nabla \btau - \nabla \bu \cdot \btau -\btau \cdot \!\!\transp{(\nabla \bu)}.$$
Most of the models of viscoelastic flows can be seen as generalizations of the Oldroyd-B model, and as we shall see, some of these generalizations admit global solutions.

\paragraph{Many objective (frame indifferent) models}

The classical way to introduce this Oldroyd-B model is to compare any elementary fluid element to a one dimensional mechanical system composed by springs and dashpots.
The derivative~$\stackrel{\triangledown}{\btau}$ is one way to extrapolate the convected derivative while preserving the invariance under galilean transformation.
There exists a one-parameter family of such models.
This parameter is usually denoted by "$a$" and the Oldroyd-B case corresponds to the case $a=1$. 

Such a models have been extensively studied.
Guillopé and Saut~\cite{Guillope-Saut4,Guillope-Saut-CRAS,Guillope-Saut3,Guillope-Saut1} proved the existence of local strong solutions.
Fernández-Cara, Guillén and Ortega~\cite{Fernandez-Guillen-Ortega-CRAS,Fernandez-Guillen-Ortega} proved local well posedness in Sobolev spaces.
In Chemin and Masmoudi~\cite{Chemin}, local and global well-posedness in critical Besov spaces were given.
In these series of papers, some global existence results hold always assuming small data.

It may be noted that in the case $a=0$ (namely the corotationel case) - and only for that case, a global existence result of weak solution has been shown, see the result of Lions and Masmoudi~\cite{Lions-Masmoudi-viscoelastique}.
 
 \paragraph{Micro-macro approach}
 
On the other hand, the Oldroyd-B model can be seen as a special case of micro-macro models.
This family of models is based on the fact that the constraint can be defined (using the formula of Kramers) from the distribution of the polymer chains.
The distribution function is itself a solution of an equation of Fokker-Planck type wherein  a spring acts.

The complete model couples the Navier-Stokes equations and this Fokker-Planck equation. A lot of local existence results are proved according to the exact form of the potential spring force, see for instance~\cite{Jourdain-Lebris-Lelievre,Masmoudi08,Renardy,Zhang-Zhang}.
Note that the Oldroyd-B model corresponds to the case where the spring force is assumed to be a linear hookean force.

Recently, Masmoudi~\cite{Masmoudi13} proved global existence of weak solutions to the FENE (Finite Extensible Nonlinear Elastic) dumbbell model.
In this model, a polymer is idealized as an elastic dumbbell consisting of two beads joined by a spring whose elongation cannot exceed a limit. The spring force therefore has a very specific shape.

\paragraph{Integral models}

Finally, the Oldroyd-B model is a special case of integral-type models.
These models are built on the natural remark expressing the fact that the fluid is a memory medium: the stress at a given time depends on all past constraints.
The Oldroyd-B model corresponds to a linear case (the influence of the Finger tensor is linear).
We show in this paper that for usual integral models, more precisely some of those including a nonlinear dependence with respect to the Finger tensor, a global existence result.\\

The situation describes before can be represented by the following picture where the Oldroyd-B model can be viewed as a particular case of some different approaches:
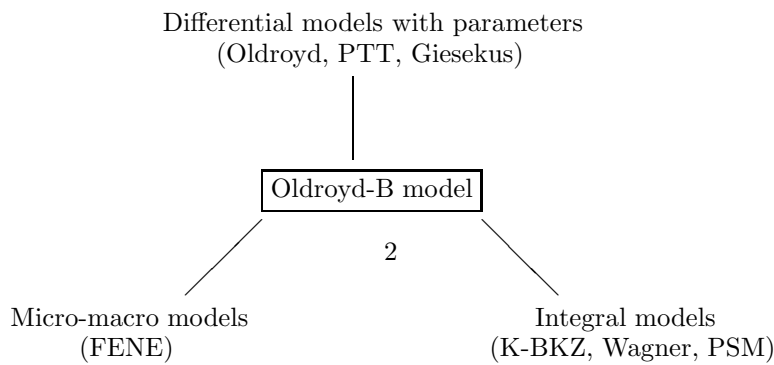
\setlength{\unitlength}{1cm}
\begin{center}
\begin{picture}(1,1)
\put(-2.5,0.5){Differential models with parameters}
\put(-1.8,0.1){(Oldroyd, PTT, Giesekus)}
\put(0,-1.2){\line(0,1){1.1}}
\put(-1.2,-1.7){\boxed{\text{Oldroyd-B model}}}
\put(-2.2,-3.0){\line(1,1){1}}
\put(-4.5,-3.4){Micro-macro models}
\put(-3.6,-3.8){(FENE)}
\put(1.7,-2.0){\line(1,-1){1}}
\put(2.4,-3.4){Integral models}
\put(1.8,-3.8){(K-BKZ, Wagner, PSM)}
\end{picture}
\end{center}
\vspace{4cm}

\paragraph{Monodimensional case - shear flows}

Some global existence results (without assuming that the data are small) already exist for some integral models: they correspond to some special flows which can be view as monodimensional cases.
Indeed, if the flow is assumed to be sufficiently simple then the Navier-Stokes equations reduce to more simple equations which automatically imply more results.
The Poiseuille flow of a KBKZ-fluid is then study in~\cite{Aarts}. They especially study the steady flow equation and its stability.
More recently, Renardy~\cite{Renardy09} proved the global existence in time of solutions to time-dependent shear flows for such integral viscoelastic behavior. The essential point in the proof is an a priori estimate for the shear stress which allows to easily deduce - in this ``one'' dimensional case - bound on the shear velocity. 

\subsection{Outline of this paper}

In the next Section (Section~\ref{part:model}), we present the model coupling the hydrodynamic Navier-Stokes equations, and the stress constitutive relation.
The Section~\ref{part:framework} is devoted to the presentation of the mathematical framework. We also give in the Section~\ref{part:framework} the main assumptions on the functional~$\mathcal F$ introduced in~\eqref{stress-integral}.
In the Section~\ref{part:result} we give the main result and describe the method for the proof.
The last three Sections (\ref{part:velocity-estimate}, \ref{part:stress-estimate} and~\ref{part:end-proof}) are devoted to the proof.
More exactly, we first give some estimates on the velocity field in the Section~\ref{part:velocity-estimate}. Next we show how to control the extra-stress using the velocity (Section~\ref{part:stress-estimate}). We finally conclude the proof in the Section~\ref{part:end-proof}.

\section{Governing equations}\label{part:model}

For a general viscoelastic and incompressible fluid, we start from the following equations for the conservation of momentum and mass
\begin{equation}\label{NS}
\begin{aligned}
& \partial_t \bu + \bu\cdot \nabla \bu + \nabla p - \eta \, \Delta \bu = \div \, \btau, \\
& \div \, \bu = 0,
\end{aligned}
\end{equation}
The two unknowns are the vector velocity field~$\bu$ and the scalar pressure~$p$.
The positive real~$\eta$ is the kinematic viscosity of the fluid.
This system is closed using a constitutive equation connecting the extra-stress~$\btau$ and the velocity gradient~$\nabla \bu$. The role of this additional contribution~$\btau$ is to take into account the past history of the fluid. It can be express by an integral with respect all past time\footnote{This case is the particular case of the separable single-integral model. We can use more general models like those given by~\eqref{stress-integral}.
In this paper, the proofs are written with a separate model~\eqref{stress} but they can easily be generalized, see the Remark~\ref{rem:non-separable-stress}.
}:
\begin{equation}\label{stress}
\btau(t,\bx) = \int_{-\infty}^{t} m(t-\sigma) \, \mathcal S\big(\bF(\sigma,t,\bx)\big)\, \mathrm d\sigma.
\end{equation}
The scalar function $m$ (the memory) and the tensorial function $\mathcal S$ are given by the properties of the fluids studied, whereas the deformation tensor $\bF$ is coupled with the velocity field of the flow. More precisely the tensor~$\bF$ satisfies the following relation
\begin{equation}\label{deformation}
\partial_t \bF + \bu \cdot \nabla \bF = \bF \cdot \nabla \bu.
\end{equation}
In this paper we are interested in the two dimensional periodical case with respect to the spatial variable:~$\bx\in \T^2$. Consequently there is no boundary condition.
We must impose the initial conditions.
For the velocity, they correspond to a given velocity at $t=0$.
For the deformation tensor~$\bF$ we give its initial value at $t=0$.
We also note that by definition of the deformation, we must have $\bF(\sigma,\sigma,\bx)=\bdelta$ for all past time~$\sigma$ and for any $\bx\in \T^2$ (the tensor~$\bdelta$ representing the identity tensor).\\

It may be more interesting to work with the new variable~$s=t-\sigma$ which represents the age instead of the parameter~$\sigma$.
We then introduce $\bG(s,t,\bx)=\bF(t-s,t,\bx)$ and the system reads
\begin{equation}\label{system1}
\left\{
\begin{aligned}
& \partial_t \bu + \bu\cdot \nabla \bu  + \nabla p - \eta \Delta \bu = \div \, \btau && \text{on $(0,+\infty)\times \T^2$,}\\
& \div \, \bu = 0 && \text{on $(0,+\infty)\times \T^2$,}\\
& \btau(t,\bx) = \int_0^{+\infty} m(s) \, \mathcal S\big(\bG(s,t,\bx)\big)\, \mathrm ds && \text{for $(t,\bx)\in(0,+\infty)\times \T^2$,}\\
& \partial_s \bG + \partial_t \bG + \bu \cdot \nabla \bG = \bG \cdot \nabla \bu && \text{on $(0,+\infty)\times (0,+\infty)\times \T^2$.}
\end{aligned}
\right.
\end{equation}
The System~\eqref{system1} is closed with the following initial conditions:
\begin{equation}\label{BC}
\bu\big|_{t=0} = \bu_0, 
\qquad
\bG\big|_{t=0} = \bG_0,
\qquad
\bG\big|_{s=0} = \bdelta.
\end{equation}

\section{Mathematical framework and assumptions}\label{part:framework}

\subsection{Tensorial analysis}

In the System~\eqref{system1}, the first equation is a vectorial equation (the velocity $\bu$ is a function with values in $\R^2$), and the two last equations are tensorial equations (the stress~$\btau$ and the deformation tensor~$\bG$ are functions with values in the set of the $2$-tensors).
In the following proofs, we need to work with the gradient of such $2$-tensors, that is with $3$-tensors, and even with $4$-tensors.
We introduce here some definitions for tensors of arbitrary order.
\begin{definition}
Let $\bA$ be a $p$-tensor and $\bB$ be a $q$-tensor.
For any $0\leq s \leq min\{p,q\}$ we define the following $(p+q-2s)$-tensor $\bA \overset{(s)}{:} \bB$ component by component:\\[-0.3cm]
\begin{equation*}
\Big( \bA \overset{(s)}{:} \bB \Big)_{i_1,...i_{p-s},j_{s+1},...,j_q} = \sum_{k_1,...,k_s}a_{i_1,...,i_{p-s},k_1,...,k_s} b_{k_1,...,k_s,j_{s+1},...,j_q}.\\[-0.2cm]
\end{equation*}
For simplicity, we will denote $\bA \overset{(0)}{:} \bB = \bA \bB$, $\bA \overset{(1)}{:} \bB = \bA \cdot \bB$ and $\bA \overset{(2)}{:} \bB = \bA:\bB$.
\end{definition}
Note also that the product $\overset{(p)}{:}$ is a scalar product on the set of the $p$-tensors. It allows us to define a generalized Froebenius norm:
\begin{definition}
The Foebenius norm of a $p$-tensor~$\bA$ is defined by $\big(\bA \overset{(p)}{:} \bA\big)^{1/2}$. It will always be denoted~$|\bA|$ (regardless of the value of $p$). Using the component of $\bA$ we have
\begin{equation*}
|\bA|^2 = \sum_{i_1,...,i_p} a_{i_1,...,i_p}^2.
\end{equation*}
\end{definition}
Finally, the Cauchy-Schwarz inequality can be easily generalized as follow:
\begin{proposition}
Let $\bA$ be a $p$-tensor, $\bB$ be a $q$-tensor and $0\leq s \leq min\{p,q\}$.
We have
\begin{equation}\label{CS}
| \bA \overset{(s)}{:} \bB | \leq |\bA| |\bB|.
\end{equation}
\end{proposition}
Note that the norms used in this proposition are not all the same: on the left hand side of the inequality~\eqref{CS}, it corresponds to the Froebenius norm on the $p+q-2s$-tensors, whereas on the right hand side member it corresponds to the Froebenius norm on the $p$-tensors and on $q$-tensors.
%

\subsection{Functional spaces}

We use the following usual notations:
\begin{itemize}
\item[$-$] For all real $s\geq 0$ and all integer $q\geq 1$, the set $W^{s,q}(\T^2)$ corresponds to the Sobolev spaces. We classically denote $L^q(\T^2)=W^{0,q}(\T^2)$ the associated Lebesgue space.
\end{itemize}
Since we will frequently use functions with values in $\R^2$ or in the space~$\mathcal L(\R^2)$ of real $2$-tensors, the usual notations will be abbreviated.
For instance, the space $(W^{1,q}(\T^2))^2$ will be denoted $W^{1,q}(\T^2)$.
Moreover, all the norms will be denoted by index, for instance like $\|\bu\|_{W^{1,q}(\T^2)}$.
\begin{itemize}
\item[$-$] The space $D^r_q(\T^2)$ stands for some fractional domain of the Stokes operator~$A_q$ in $L^q(\T^2)$ (cf. Section 2.3 in~\cite{Danchin}). Its norm is defined by
$$
\|\bv\|_{D^r_q(\T^2)} := \|\bv\|_{L^q(\T^2)} + \Big( \int_0^{+\infty} \|A_q \mathrm e^{-tA_q}\bv \|_{L^q(\T^2)}^{r} \, \mathrm dt \Big)^{1/r}.
$$
\end{itemize}
Roughly, the vector-fields of~$D^r_q(\T^2)$ are vectors which have $2-\frac{2}{r}$ derivatives in $L^q(\T^2)$ and are divergence-free.
It may be identified with Besov spaces. It also can be view as an interpolate space between $L^q(\T^2)$ and the domain of the Stokes operator~$D(A_q)$, see~\cite{Danchin}.
\begin{itemize}
\item[$-$] The notation of kind $L^r(0,T;W^{1,q}(\T^2))$ denotes the space of $r$-integrable functions on $(0,T)$, with values in~$W^{1,q}(\T^2)$.
Similarly, expression like $g\in L^\infty(\R^+;L^r(0,T;L^q(\T^2)))$ means that
$$
\sup_{s\in \R^+} \Big( \int_0^T \|g(s,t,\cdot)\|_{L^q(\T^2)}^r \, \mathrm dt \Big)^{\frac{1}{r}} < +\infty.
$$
\item[$-$] Finally let us write~$\mathfrak P$ for the orthogonal projector in~$L^2(\T^2)$ onto the set of the divergence-free vectors fields of~$L^2(\T^2)$.
\end{itemize}

\subsection{Assumptions}\label{sec:assumption}

In this Section we present the assumptions that we need for the proof. These assumptions concern the functions~$m$ and~$\mathcal S$ introduced in the extra-stress expression~\eqref{stress}:
\begin{enumerate}
\item[(H1)] $m:s\in \R^+\longmapsto m(s)\in \R$ is measurable, decreasing, positive and satisfies $\int_0^\infty m(s)\, \mathrm ds = 1$;
\item[(H2)] $\mathcal S:\bG\in \mathcal L(\R^2) \longmapsto \mathcal S(\bG) \in \mathcal L(\R^2)$ is of class $\mathcal C^1$ and satisfies
\begin{itemize}
\item[--] There exists $\mathcal S_\infty\geq 0$ such that for all $\bG\in \mathcal L(\R^2)$ we have $\dsp |\mathcal S(\bG)|\leq \mathcal S_\infty$;
\item[--] There exists $\mathcal S_\infty'\geq 0$ such that for all $\bG\in \mathcal L(\R^2)$ we have $ |\bG| |\mathcal S'(\bG)|\leq \mathcal S_\infty'$.
\end{itemize}
\end{enumerate}
As precised above the matricial norms used here correspond to the Froebenius norms. We take care of the fact that the derivative~$\mathcal S'(\bG)$ may be represented by a tensor of order~$4$: $(S'(\bG))_{ijk\ell}$ corresponds to the derivation of $(\mathcal S(\bG))_{k\ell}$ with respect to the component $\bG_{ij}$.

\paragraph{Notes on the assumptions}
\begin{itemize}
\item[$\checkmark$] The first assumption (H1) is related to the memory function~$m$. It is linked to the principle of fading memory, see~\cite{Coleman}. Usually, the memory function is a combination of exponentially decreasing functions which satisfies the assmption (H1).
Note that in some cases, the memory function is described as an infinite sum of exponentially decreasing functions. This is the case of the Doi-Edwards model, see~\cite{Doi}. Despite the singularity that has such function at~$0$, it satisfies the hypothesis (H1).

\item[$\checkmark$] The second assumption concerns the function~$\mathcal S$. It is satisfied by a lot of classical integral models.
For example, in the two dimensional case, usual models read
$$\mathcal S(\bG) = h(I_1)\transp{\bG}\cdot \bG,$$
where $I_1=\mathrm{Tr}(\transp{\bG}\cdot \bG)$ is the only invariant which is of interest (the other one is given by~$\mathrm{det}(\transp{\bG}\cdot \bG)$. It is equal to~$1$ since the flow is assumed to be incompressible).
For such a cases, the assumption (H2) is equivalent to
\begin{itemize}
\item[--] There exists $C\geq 0$ such that for all $x\geq 0$ we have $x|h(x)| \leq C$;
\item[--] There exists $C'\geq 0$ such that for all $x\geq 0$ we have $x^2|h'(x)|\leq C'$.
\end{itemize}
For the PSM model ($h(x) \approx \frac{1}{1+x}$) or for the Wagner model ($h(x) \approx \mathrm e^{-\sqrt{x}}$) these assumptions are clearly verified.
\end{itemize}

Nevertheless, we note that the Oldroyd-B model, which corresponds to the case $m(s)=\mathrm e^{-s}$ and $\mathcal S(\bG) = \transp{\bG}\cdot \bG - \bdelta$, does not satisfied the assumption~(H2).
The study presented here does not cover such Oldroyd models: the global result in this case remaining an open question. 

\begin{remark}\label{rem:non-separable-stress}
If we want to use a non-separable integral law like
\begin{equation}
\btau(t,\bx) = \int_0^{+\infty} \mathcal F\big(s,\bG(s,t,\bx)\big)\, \mathrm ds,
\end{equation}
the assumptions (H1) and (H2) become
\begin{enumerate}
\item[--] There exists $m_1\in L^1(\R^+)$ such that for $(s,\bG)\in \R^+\!\times\! \mathcal L(\R^2)$ we have $\dsp |\mathcal F\big(s,\bG \big)| \leq m_1(s)$;
\item[--] There exists $m_2\in L^1(\R^+)$ decreasing such that for $(s,\bG)\in \R^+\!\times\! \mathcal L(\R^2)$ we have $\dsp |\bG| |\partial_{\bG}\mathcal F\big(s,\bG \big)| \leq m_2(s)$.
\end{enumerate}
All the proofs remain unchanged.
\end{remark}

\section{Main result}\label{part:result}

\begin{theorem}\label{theorem}
Let $q$ and~$r$ be two integers such that $\frac{1}{q}+\frac{1}{r}<\frac{1}{2}$.
We assume that the initial conditions~$\bu_0$ and~$\bG_0$ satisfy
$$
\bu_0\in D_q^r(\T^2), \qquad \bG_0\in L^\infty(\R^+;W^{1,q}(\T^2)) \cap W^{1,\infty}(\R^+;L^q(\T^2)),
$$
and there exists~$\mu>0$ such that $\det \bG_0\geq \mu$ on~$\R^+\times \T^2$.
Let $\eta>0$, $m$ satisfying~(H1), $\mathcal S$ satisfying~(H2) and $T>0$ be arbitrary.

There exists a constant $C$ depending only on the norm of the initial data, $q$, $r$, $\mu$, $\eta$, $\mathcal S_\infty$, $\mathcal S_\infty'$ and~$T$ with $C$ bounded for bounded~$T$, and a unique solution $(\bu,p,\btau, \bG)$ of \eqref{system1}--\eqref{BC} such that
\begin{equation*}
\begin{aligned}
& \|\nabla^2 \bu\|_{L^r(0,T;L^q(\T^2))} \leq C,
&& \|\nabla \bu\|_{L^\infty((0,T)\times \T^2)} \leq C,\\
& \|\nabla \btau\|_{L^r(0,T;L^q(\T^2))} \leq C,
&& \|\btau\|_{L^\infty((0,T)\times \T^2)} \leq C,
\end{aligned}
\end{equation*}
\vspace{-0.2cm}
\begin{equation*}
\text{and} \quad \int_0^T \int_0^\infty m(s) \Big\| \frac{\nabla \bG}{|\bG|} \Big\|_{L^q(\T^2)}^r(s,t) \, \mathrm ds \mathrm dt \leq C,
\end{equation*}
hold.
\end{theorem}

\begin{remark}
\begin{itemize}
\item[]
\item[$\checkmark$] The pressure~$p$ is a Lagrange multiplier associated to the divergence free constraint.  
It can be solve using the Riesz transforms. More precisely, taking the divergence of the first equation of the System~\eqref{system1} we use the periodic boundary conditions to have
\begin{equation}\label{pressure}
p=-(-\Delta)^{-1}\div\div\, (\btau-\bu\otimes\bu).
\end{equation}
From the Theorem~\ref{theorem}, the solutions of the System~\eqref{system1} discussed in this paper have $\btau-\bu \otimes \bu$ in $L^\infty(0,T;L^2(\T^2))$. The pressure in the solution of~\eqref{system1} is meant to be given by~\eqref{pressure}.
\item[$\checkmark$] In many application, the fluid is assumed to be initially quiescent. In that case, we have $\bG_0=\bdelta$ and $\det\bG_0=1$. Moreover, we will see that the quantity $\det\bG$ is only convected by the flow. If the fluid is assumed to be at rest in the past (that is for $s$ large enough), then we always have  $\det\bG_0=1$. The assumption on the positivity of~$\det\bG_0$ allows us consider, for instance, such cases.
\end{itemize}
\end{remark}

In the following, we will denote by~$C$ constants that may depend on the initial conditions, on the viscosity~$\eta$, on some integer $r$, $q$, on the bounds~$\mathcal S_\infty$ and~$\mathcal S_\infty'$, on the constant~$\mu$, and on the time~$T$.
These constants will always be bounded for bounded~$T$.

\paragraph{Sketch of the proof}

Using the assumptions given in the Theorem~\ref{theorem}, the local existence is proved in~\cite{Chupin3}. It is based on a fixed point argument and some estimates. The existence time is small since we need some contraction in the fixed point Theorem.
Nevertheless, to obtain the local existence we do not need assumption~(H2): we only assume that the function~$\mathcal S$ is of class~$\mathcal C^1$.

The purpose of this article is to establish additional bounds using the additional assumption~(H2).
We then consider a solution~$(\bu,p,\btau,\bG)$ to the System~\eqref{system1}--\eqref{BC} in $[0,T]$ with the regularity proved in~\cite{Chupin3}:
\begin{equation*}
\begin{array}{ll}
 \bu \in L^r(0,T;W^{2,q}(\T^2)),
& \partial_t \bu \in L^r(0,T;L^q(\T^2)),\\
 \btau \in L^\infty(0,T;W^{1,q}(\T^2)),
& \partial_t \btau \in L^r(0,T;L^q(\T^2)),\\
 \bG \in L^\infty(\R^+ \! \times \! (0,T);W^{1,q}(\T^2)), 
& \partial_s \bG, \, \partial_t \bG \in L^\infty(\R^+;L^r(0,T;L^q(\T^2))).
\end{array}
\end{equation*}
The next steps are to obtain estimates on this solution.

Roughly speaking the first part of the assumption~(H2) implies that the extra-stress~$\btau$ is $L^\infty$-bounded. The second part of the assumption~(H2) gives a control of $\nabla \btau$ with respect to $\frac{\nabla \bG}{\bG}$.
These control on the extra-stress will be transform into controls on the velocity using the Navier-Stokes equations.
Finally the equation on~$\bG$ allows us to deduce a bound on~$\frac{\nabla \bG}{\bG}$.
 
\section{A priori estimates for the spatial gradient of the velocity}\label{part:velocity-estimate}

The following key result is a direct consequence of the assumptions (H1) and (H2) on the stress tensor by means of the function~$\mathcal S$:
\begin{lemma}\label{lemma1}
We have the following $L^\infty$ bound:
\begin{equation}\label{estimate1}
\|\btau\|_{L^\infty((0,T)\times \T^2)} \leq \mathcal S_\infty.
\end{equation}
\end{lemma}
We can prove that the velocity field is also bounded:
\begin{lemma}\label{lemma2}
There exists a constant $C$ such that
\begin{equation}\label{estimate2}
\|\bu\|_{L^\infty((0,T)\times \T^2)} \leq C.
\end{equation}
\end{lemma}

\proof
In the one hand, we use the local in time result to obtain a bound for~$\bu$ in $(0,T_0)\times \T^2$ for some~$T_0>0$.
Indeed the local existence result gives a bound for $\bu$ in $L^r(0,T_0;W^{2,q}(\T^2))$ and a bound for $\partial_t\bu$ in $L^r(0,T_0;L^q(\T^2))$. For $r\geq 2$, by a Aubin-Simon Theorem (see~\cite{Simon}) this implies a bound for $\bu$ in $\mathcal C(0,T_0;W^{1,q}(\T^2))$ which, for $q>2$, provides the $L^\infty((0,T_0)\times\T^2)$ bound on~$\bu$.
\par
In the other hand, a result proved by Constantin and Seregin (see~\cite[Prop. 2.4]{Constantin-Seregin}) gives a $L^\infty$-bound for the solution~$\bu$ to the Navier-Stokes Equations~\eqref{NS} in $(\sigma,T)\times \T^2$ for any~$\sigma>0$, as soon as $\btau$ is bounded in some $L^4((0,T)\times \T^2)$.
\par
Taking $\sigma=T_0/2$ this allows to conclude the proof of the lemma~\ref{lemma2}.
\cqfd

\begin{lemma}\label{lemma3}
For all $1<q,r<+\infty$ there exists a constant $C$ such that
\begin{equation}\label{estimate3}
\|\nabla \bu\|_{L^r(0,T;L^q(\T^2))} \leq C.
\end{equation}
\end{lemma}

\proof
The proof is based on the integral representation of the solution to the Navier-Stokes Equation~\eqref{NS}:
\begin{equation}\label{NSInt}
\nabla \bu(t,\bx) =\mathrm e^{\eta t \Delta} \nabla \bu_0 + \int_0^t \mathrm e^{\eta (t-\sigma) \Delta} \mathfrak P \Delta(\btau - \bu\otimes \bu)(\sigma,\bx) \, \mathrm d\sigma.
\end{equation}
We use the fact that the linear operator
$\mathfrak T:\bbf \longmapsto \int_0^t \mathrm e^{\eta (t-\sigma) \Delta} \Delta \bbf(\sigma) \, \mathrm d\sigma$,
is bounded in $L^r(0,T;L^q(\T^2))$ for $1<q,r<+\infty$, see~\cite[p. 64]{Lemarie}.
The previous lemmas~\ref{lemma1} and~\ref{lemma2} give estimates for $\bbf=\mathfrak P (\btau - \bu\otimes \bu)$ in $L^r(0,T;L^q(\T^2))$ for any $1<q,r<+\infty$, that complete the proof of this lemma~\ref{lemma3}.
\cqfd

\begin{proposition}\label{proposition1}
For $\frac{1}{q}+\frac{1}{r}<\frac{1}{2}$ there exists a constant $C$ such that for all $t\in (0,T)$
\begin{equation}\label{estimate-2}
\|\nabla \bu(t,\cdot)\|_{L^\infty(\T^2)} \leq C + C \ln(\mathrm e+\|\nabla \btau\|_{L^r(0,t;L^q(\T^2))}),
\end{equation}
\begin{equation}\label{estimate-3}
\|\nabla^2 \bu\|_{L^r(0,t;L^q(\T^2))} \leq C + C\|\nabla \btau\|_{L^r(0,t;L^q(\T^2))}.
\end{equation}
\end{proposition}

\proof
The proof is also based on the integral representation~\eqref{NSInt}.
We will use the following result about the kernel of the heat equation (see~\cite{Lemarie}):
Fisrt of all, if $\bbf(\sigma,\cdot)\in L^\infty(\T^2)$ for all $\sigma\in (0,T)$ then we have, for all $\sigma\in (0,T)$:
\begin{equation}\label{ker1}
\|\mathrm e^{\eta (t-\sigma) \Delta} \Delta \bbf(\sigma,\cdot) \|_{L^{\infty}(\T^2)} \leq C (t-\sigma)^{-1} \|\bbf(\sigma,\cdot) \|_{L^{\infty}(\T^2)}.
\end{equation}
Next, if $\bbf(\sigma,\cdot)\in L^q(\T^2)$ for all $\sigma\in (0,T)$ and $1<q<\infty???$ then we have, for all $\sigma\in (0,T)$:
\begin{equation}\label{ker2}
\|\mathrm e^{\eta (t-\sigma) \Delta} \Delta \bbf(\sigma,\cdot) \|_{L^{\infty}(\T^2)} \leq C (t-\sigma)^{-\frac{q+2}{2q}} \|\nabla \bbf(\sigma,\cdot) \|_{L^q(\T^2)}.
\end{equation}
Denoting $\bbf =  \mathfrak P(\btau - \bu\otimes \bu)$, the expression~\eqref{NSInt} reads, for any $0<t^\star<t$:
\begin{equation}\label{eq:NSgradInt}
\begin{aligned}
\nabla \bu(t,\bx) =\mathrm e^{\eta t \Delta} \nabla \bu_0 + \int_0^{t-t^\star} \mathrm e^{\eta (t-\sigma) \Delta} \Delta \bbf (\sigma,\bx) \, \mathrm d\sigma + \int_{t-t^\star}^t \mathrm e^{\eta (t-\sigma) \Delta} \Delta \bbf (\sigma,\bx) \, \mathrm d\sigma.
\end{aligned}
\end{equation}
We take the $L^\infty$-norm with respect to the spatial variable and we use the results~\eqref{ker1} and~\eqref{ker2} to obtain
\begin{equation}\label{eq:NSgradInt}
\begin{aligned}
\|\nabla \bu(t,\cdot)\|_{L^\infty(\T^2)} \leq C + C \int_0^{t-t^\star} (t-\sigma)^{-1} \| \bbf (\sigma,\cdot) \|_{L^\infty(\T^2)} \, \mathrm d\sigma
+ C \int_{t-t^\star}^t (t-\sigma)^{-\frac{q+2}{2q}} \| \nabla \bbf(\sigma,\cdot) \|_{L^q(\T^2)} \, \mathrm d\sigma.
\end{aligned}
\end{equation}
Using the Hölder inequality, we deduce
\begin{equation}\label{eq:36}
\begin{aligned}
\|\nabla \bu(t,\cdot)\|_{L^\infty(\T^2)}
& \leq  C + C \ln\Big(\frac{t}{t^\star}\Big) \| \bbf \|_{L^\infty((0,T)\times \T^2)}
+ C \Big( \int_{t-t^\star}^t (t-\sigma)^{-\frac{q+2}{2q}\frac{r}{r-1}}  \, \mathrm d\sigma \Big)^{\frac{r-1}{r}} \| \nabla \bbf\|_{L^r((0,t);L^q(\T^2))}\\
& \leq  C + C \ln\Big(\frac{t}{t^\star}\Big) \| \bbf \|_{L^\infty((0,T)\times \T^2)}
+ C {t^\star}^{\alpha} \| \nabla \bbf \|_{L^r((0,t);L^q(\T^2))}.
\end{aligned}
\end{equation}
where $\alpha=\frac{1}{2}-\frac{1}{q}-\frac{1}{r}$ is positif due to the assumption $\frac{1}{q}+\frac{1}{r} <\frac{1}{2}$.
According to the lemmas~\ref{lemma1} and~\ref{lemma2} we know that $\|\bbf\|_{L^\infty((0,T)\times \T^2)} \leq C$. In the same way, according to the lemmas~\ref{lemma1}, \ref{lemma2} and~\ref{lemma3}, we have $\|\nabla \bbf\|_{L^r(0,t;L^q(\T^2))} \leq C+ C\|\nabla \btau\|_{L^r(0,t;L^q(\T^2))}$.
The inequality~\eqref{eq:36} reads
\begin{equation}\label{eq:37}
\begin{aligned}
\|\nabla \bu(t,\cdot)\|_{L^\infty(\T^2)}
& \leq  C + C \ln\Big(\frac{t}{t^\star}\Big)
+ C {t^\star}^{\alpha} \| \nabla \btau \|_{L^r((0,t);L^q(\T^2))}.
\end{aligned}
\end{equation}
We now choose
$$t^\star = \min\{\mathrm e^{-1}, \| \nabla \btau \|_{L^r((0,t);L^q(\T^2))}^{-1/\alpha}\} \, t.$$
Since $\mathrm e^{-1}<1$ we have $0<t^\star<t$, and since $\alpha<\frac{1}{2}$ we have $\ln\big(\frac{t}{t^\star}\big) \leq \frac{1}{\alpha}\ln\big( \mathrm e + \| \nabla \btau \|_{L^r((0,t);L^q(\T^2))}\big)$. The estimate~\eqref{eq:36} gives the first result~\eqref{estimate-2} of the proposition~\ref{proposition1}.
\par
To prove the second inequality~\eqref{estimate-3} of the proposition~\ref{proposition1}, we take the spatial gradient of the expression~\eqref{NSInt}:
\begin{equation}\label{eq:NSlapInt}
\nabla^2 \bu(t,\bx) =\mathrm e^{\eta t \Delta} \nabla^2 \bu_0 + \int_0^t \mathrm e^{\eta (t-\sigma) \Delta} \Delta \nabla \bbf(\sigma,\bx) \, \mathrm d\sigma.
\end{equation}
Taking the $L^r(0,t;L^q(\T^2))$ norm, the initial term $\mathrm e^{\eta t \Delta} \nabla^2 \bu_0$ exactly corresponds to the norm of $\bu_0$ in the space~$D_q^r(\T^2)$. The integral term is controled using the boundness of the operateur~$\mathfrak T$ introduced in the proof of the lemma~\ref{lemma3}.
We note once again the control of $\|\nabla \bbf\|_{L^r(0,t;L^q(\T^2))}$ using~$\|\nabla \btau\|_{L^r(0,t;L^q(\T^2))}$.
\cqfd

\section{Control of the stress gradient}\label{part:stress-estimate}

\begin{lemma}\label{lemma7}
There exists a constant $C$ such that for all $(s,t,\bx)\in \R^+\times (0,T)\times \T^2$ we have
\begin{equation}\label{estimate7}
|\bG(s,t,\bx)| \geq C >0.
\end{equation}
\end{lemma}

\proof
By assumption, for all $(s,\bx)\in (0,+\infty)\times \T^2$ we have:
\begin{equation}\label{eq:12}
\det(\bG(s,0,\bx))\geq \mu>0.
\end{equation}
Moreover we have $\bG|_{s=0}=\bdelta$ so that for all $(t,\bx)\in (0,T)\times \T^2$ we have:
\begin{equation}\label{eq:13}
\det(\bG(0,t,\bx)) =1.
\end{equation}
\par
A simple calculation shows that the quantity~$\det(\bG)$ satisfies
$$\mathcal D \det(\bG) =  \!\div\, \bu \, \det(\bG)= 0,$$
where~$\mathcal D$ refers to the one order derivating operator~$\mathcal D=\partial_s+\partial_t+\bu\cdot \nabla$.
The value~$\det(\bG)$ is then constant along the characteristic lines.
Since all the characteristic lines start from the lines $\{s=0\}$ and $\{t=0\}$ we deduce from~\eqref{eq:12} and~\eqref{eq:13} that $\det(\bG)\geq \min(\mu,1)$ on $\R^+\times (0,T)\times \T^2$.
\par
Due to the inequality of arithmetic and geometric means, we have 
$$|\bG|^2 = \mathrm{Tr}(\transp{\bG}\cdot\bG) \geq 2\sqrt{\det(\transp{\bG}\cdot\bG)} = 2|\det(\bG)| \geq 2\min(\mu,1),$$
that concludes the proof of the lemma~\ref{lemma7}.
\cqfd

Since $\nabla \bG \in L^\infty(\R^+\!\times\! (0,T);L^q(\T^2))$, this lemma~\ref{lemma7} implies $\dsp \frac{\nabla \bG}{|\bG|} \in L^\infty(\R^+\!\times\! (0,T);L^q(\T^2))$.
We use this quantity to estimate the gradient of the stress:
\begin{lemma}\label{lemma62}
For $1<q,r<+\infty$, and for all $t\in (0,T)$ we have
\begin{equation}\label{estimate7}
\|\nabla \btau \|_{L^r(0,t;L^q(\T^2))}^r
\leq
\mathcal S_\infty' \int_0^t \int_0^\infty m(s) \Big\| \frac{\nabla \bG}{|\bG|} \Big\|_{L^q(\T^2)}^r(s,t) \, \mathrm ds \mathrm dt.
\end{equation}
\end{lemma}

\proof
To obtain estimate~\eqref{estimate7}, we first derivate the stress tensor~$\btau$ with respect to the spatial coordinates:
\begin{equation}
\nabla \btau(t,\bx) = \int_0^\infty m(s) \, \mathcal S'(\bG(s,t,\bx)) : \nabla \bG(s,t,\bx)\, \mathrm ds.
\end{equation}
Using the assumption~(H2), we write
\begin{equation}
|\nabla \btau(t,\bx)| \leq \mathcal S_\infty' \int_0^\infty m(s) \,  \Big|\frac{\nabla \bG(s,t,\bx)}{|\bG(s,t,\bx)|}\Big| \, \mathrm ds.
\end{equation}
From the triangular inequality we deduce that
\begin{equation}
\begin{aligned}
\|\nabla \btau \|_{L^r(0,t;L^q(\T^2))}
& \leq \mathcal S_\infty' \int_0^\infty m(s) \,  \Big\|\frac{\nabla \bG}{|\bG|}\Big\|_{L^r(0,t;L^q(\T^2))}(s) \, \mathrm ds \\
& \leq \mathcal S_\infty' \int_0^\infty m(s) \bigg( \int_0^t \Big\|\frac{\nabla \bG}{|\bG|}\Big\|_{L^q(\T^2)}^r(s,\sigma) \, \mathrm d\sigma \bigg)^{1/r} \mathrm ds.
\end{aligned}
\end{equation}
Writing $m(s)=m(s)^{1-\frac{1}{r}} \times m(s)^{\frac{1}{r}}$ we apply the Hölder inequality to deduce the estimate~\eqref{estimate7} and conclude the proof of the lemma~\ref{lemma62}.
\cqfd
It is then natural to define, for all time $t\in (0,T)$ the value
\begin{equation}\label{def:y}
y(t) = \int_0^t \int_0^\infty m(s) \Big\| \frac{\nabla \bG}{|\bG|} \Big\|_{L^q(\T^2)}^r(s,\sigma) \, \mathrm ds \mathrm d\sigma.
\end{equation}
The following lemma gives an differential inequation about this quantity:
\begin{proposition}\label{proposition2}
For all integers $q$, $r$ such that $\frac{1}{q}+\frac{1}{r}<\frac{1}{2}$ and for all $t\in [0,T]$ the quantity~$y(t)$ introduced by~\eqref{def:y} satisfies
\begin{equation}
y'(t) \leq C + y(t) + C y(t) \|\nabla \bu \|_{L^\infty((0,T)\times \T^2)} + C \|\nabla^2 \bu \|_{L^r(0,t;L^q(\T^2))}^r.
\end{equation}
\end{proposition}

\proof
The Equation satisfied by~$\bG$ reads
\begin{equation}\label{eq:G}
\mathcal D \bG = \bG \cdot \nabla \bu
\qquad
\text{on $(0,+\infty)\times (0,T)\times \T^2$,}
\end{equation}
where we recall that~$\mathcal D$ corresponds to the operator~$\mathcal D=\partial_s+\partial_t+\bu\cdot \nabla$. We take the scalar product of the Equation~\eqref{eq:G} by $-q|\nabla \bG|^{q}|\bG|^{-q-2}\bG$:
\begin{equation}
|\nabla \bG|^{q} \mathcal D |\bG|^{-q} = -q|\nabla \bG|^{q}|\bG|^{-q-2} (\bG \cdot \nabla \bu):\bG.
\end{equation}
Using the generalised Cauchy-Schwarz inegality~\eqref{CS}, we deduce
\begin{equation}\label{eq:10}
|\nabla \bG|^{q} \mathcal D |\bG|^{-q} \leq q|\nabla \bG|^{q}|\bG|^{-q} |\nabla \bu|.
\end{equation}
Next we take the spatial derivative of the Equation~\eqref{eq:G}.
We obtain the following $3$-tensor equation
\begin{equation}
\mathcal D \nabla \bG = \nabla \bG \cdot \nabla \bu + (\bG\cdot \nabla^2\bu)^\dag - \nabla \bu \cdot \nabla \bG.
\end{equation}
More precisely, the component $(i,j,k)$ of this equation reads 
\begin{equation}
\mathcal D \partial_i \bG_{jk} = \partial_i\bG_{j\ell} \partial_\ell \bu_k + \bG_{j\ell} \partial_\ell \partial_i \bu_k - \partial_i\bu_\ell\partial_\ell \bG_{jk}.
\end{equation}
Taking the scalar product of this equation by $q|\bG|^{-q}|\nabla \bG|^{q-2} \nabla \bG$ and using the Cauchy-Schwarz inequality we deduce
\begin{equation}\label{eq:20}
|\bG|^{-q} \mathcal D |\nabla \bG|^q \leq 2q|\nabla \bG|^{q}|\bG|^{-q}|\nabla \bu| + q |\nabla \bG|^{q-1} |\bG|^{-(q-1)} |\nabla^2\bu|.
\end{equation}
Adding these inequality~\eqref{eq:20} with the inequality~\eqref{eq:10} we deduce
\begin{equation*}
\mathcal D \big(|\nabla \bG|^q|\bG|^{-q}\big) \leq 3q|\nabla \bG|^{q}|\bG|^{-q}|\nabla \bu| + q |\nabla \bG|^{q-1} |\bG|^{-(q-1)} |\nabla^2\bu|.
\end{equation*}
Integrating with respect to the spatial variable we obtain
\begin{equation*}
\partial_s \Big\|\frac{\nabla \bG}{|\bG|} \Big\|_{L^q(\T^2)}^q + \partial_t \Big\|\frac{\nabla \bG}{|\bG|} \Big\|_{L^q(\T^2)}^q
\leq
3q \int_{\T^2} \Big|\frac{\nabla \bG}{|\bG|} \Big|^q |\nabla \bu| + q\int_{\T^2} \Big|\frac{|\nabla \bG|}{|\bG|}\Big|^{q-1} |\nabla^2\bu|.
\end{equation*}
We now use the Hölder inequality to write
\begin{equation}\label{estimate5}
\partial_s \Big\|\frac{\nabla \bG}{|\bG|} \Big\|_{L^q(\T^2)}^q + \partial_t \Big\|\frac{\nabla \bG}{|\bG|} \Big\|_{L^q(\T^2)}^q
\leq
3q \Big\|\frac{\nabla \bG}{|\bG|} \Big\|_{L^q(\T^2)}^q \|\nabla \bu\|_{L^\infty(\T^2)} + q \Big\|\frac{|\nabla \bG|}{|\bG|}\Big\|_{L^q(\T^2)}^{q-1} \|\nabla^2\bu\|_{L^q(\T^2)}.
\end{equation}
We multiply~\eqref{estimate5} by $\displaystyle \frac{r}{q}\Big\|\frac{\nabla \bG}{|\bG|}\Big\|_{L^q(\T^2)}^{r-q}$ to have
\begin{equation}
\partial_s \Big\|\frac{\nabla \bG}{|\bG|}\Big\|_{L^q(\T^2)}^r + \partial_t \Big\|\frac{\nabla \bG}{|\bG|}\Big\|_{L^q(\T^2)}^r
\leq 3r \Big\|\frac{\nabla \bG}{|\bG|}\Big\|_{L^q(\T^2)}^r \|\nabla \bu \|_{L^\infty(\T^2)} + r \Big\|\frac{\nabla \bG}{|\bG|}\Big\|_{L^q(\T^2)}^{r-1} \|\nabla^2 \bu \|_{L^q(\T^2)}.
\end{equation}
Using the Young inequality we obtain
\begin{equation}
\partial_s \Big\|\frac{\nabla \bG}{|\bG|}\Big\|_{L^q(\T^2)}^r + \partial_t \Big\|\frac{\nabla \bG}{|\bG|}\Big\|_{L^q(\T^2)}^r
\leq 3r \Big\|\frac{\nabla \bG}{|\bG|}\Big\|_{L^q(\T^2)}^r \|\nabla \bu \|_{L^\infty(\T^2)} + \Big\|\frac{\nabla \bG}{|\bG|}\Big\|_{L^q(\T^2)}^r + (r-1)^{r-1} \|\nabla^2 \bu \|_{L^q(\T^2)}^r.
\end{equation}
We multiply by~$m(s)$ and integrate for $s\in (0,+\infty)$. Assuming (H1) we deduce that the first term is non negative (we also recall that $\bG\big|_{s=0}=\bdelta$) and we obtain
\begin{equation}
y'' \leq 3r\, y' \|\nabla \bu \|_{L^\infty(\T^2)} + y' + (r-1)^{r-1} \|\nabla^2 \bu \|_{L^q(\T^2)}^r.
\end{equation}
Integrating now with respect to time in $(0,t)$, with $0\leq t\leq T$ we deduce
\begin{equation}\label{eq:123}
y'(t) \leq 3r\, y(t) \|\nabla \bu \|_{L^\infty((0,T)\times \T^2)} + y'(0) + y(t) + (r-1)^{r-1} \|\nabla^2 \bu \|_{L^r(0,t;L^q(\T^2))}^r.
\end{equation}
The value of $y'(0)$ is given with respect to the initial condition~$\bG_{\text{old}}$:
$$
y'(0) = \int_0^\infty m(s) \Big\| \frac{\nabla \bG_0}{|\bG_0|} \Big\|_{L^q(\T^2)}^r(s) \, \mathrm ds.
$$
We will note that $y'(0)$ is bounded since $\bG_0\in L^\infty(\R^+;W^{1,q}(\T^2))$ and $|\bG_0|\geq\sqrt{2\min(\mu,1)}$ on~$\T^2$:
$$
y'(0) \leq \frac{1}{(2\min(\mu,1))^{r/2}} \|\bG_0\|_{L^\infty(\R^+;W^{1,q}(\T^2))}.
$$
The estimate~\eqref{eq:123} takes the form required in the proposition~\ref{proposition2}
\cqfd

\section{Conclusion: proof of the Theorem~\ref{theorem}}\label{part:end-proof}

\begin{proposition}
The function~$y$ introduced by~\eqref{def:y} satisfies the following inequation on $(0,T)$:
\begin{equation}\label{eq:789}
y' \leq C(\mathrm e+y)\ln(\mathrm e+y).
\end{equation}
That implies $y \leq C$ on $(0,T)$.
\end{proposition}

\proof
In terms of function~$y$, the lemma~\ref{lemma62} reads: for all $t\in (0,T)$ we have
\begin{equation}
\|\nabla \btau \|_{L^r(0,t;L^q(\T^2))}^r \leq C y(t).
\end{equation}
Consequently, the estimates~\eqref{estimate-2} and~\eqref{estimate-3} of the proposition~\ref{proposition1} can be written as 
\begin{equation}\label{estimate-2b}
\|\nabla \bu(t,\cdot)\|_{L^\infty(\T^2)} \leq C + C \ln(\mathrm e+y(t)),
\end{equation}
\begin{equation}\label{estimate-3b}
\|\nabla^2 \bu\|_{L^r(0,t;L^q(\T^2))}^r \leq C + Cy(t).
\end{equation}
Using the proposition~\ref{proposition2} we deduce that the function~$y$ satisfies the following inequality on $(0,T)$:
$$
y' \leq C+Cy+Cy\ln(\mathrm e+y),
$$
that we can rewrite, up to a change of constants~$C$, as~\eqref{eq:789}.
\par
Since all the solutions of this Equation~\eqref{eq:789} are bounded for finite time:
$$y(t)\leq \mathrm e^{\mathrm e^{Ct}} \qquad \text{for all $t\in (0,T)$},$$
the proof is complete.
\cqfd


\bibliographystyle{plain}
\bibliography{biblio-global}

\end{document}